\RequirePackage{fix-cm}

\documentclass[smallextended]{svjour3}       
\smartqed  
\usepackage{graphicx}

\usepackage{epsfig, graphicx, xcolor}
\usepackage{amssymb}
\usepackage{amsmath, amsfonts}

\newtheorem {cor} {Corollary}

\newcommand {\ttbox} [1] {\mbox {\ttfamily*1}} 

\begin{document}

	\title{Inequalities of Hardy-Littlewood-Polya type for functions of operators and their applications}
	\author{Vladyslav Babenko \and Yuliya Babenko \and Nadiia Kriachko}
	
	\institute{V. Babenko \at
              Department of Mathematics and Mechanics, Dnepropetrovsk National University, Gagarina pr., 72, Dnepropetrovsk, 49010, UKRAINE\\
              \email{babenko.vladislav@gmail.com}           
           \and
           Y. Babenko \at
              Department of Mathematics, Kennesaw State University, 1100 South Marietta Pkwy, MD \# 9085, Marietta, GA, 30060, USA\\
                \email{ybabenko@kennesaw.edu}
              \and
              N. Kriachko \at
              Department of Mathematics and Mechanics, Dnepropetrovsk National University, Gagarina pr., 72, Dnepropetrovsk, 49010, UKRAINE\\
               \email{nadiakriachko@gmail.com}
}

\date{Received: date / Accepted: date}

\titlerunning{Inequalities for functions of operators}

			\maketitle
	\begin{abstract} In this paper, we derive a generalized multiplicative Hardy-Littlewood-Polya type inequality, as well as several related additive inequalities, for functions of operators in Hilbert spaces. In addition, we find the modulus of continuity of a function of an operator on a class of elements defined with the help of another function of the operator.
	We then apply the results to solve the following problems: (i) the problem of approximating a function of an unbounded self-adjoint operator by bounded operators, (ii) the problem of best approximation of a certain class of elements from a Hilbert space by another class,
 and (iii) the problem of optimal recovery of an operator on a class of elements given with an error.



	\end{abstract}

\section{Definitions, notation, and statements of main problems}\label{S1}
 Let $X,\, Y$ be Banach spaces and $A\,:\, X\to Y$ be an operator (not necessarily linear) with domain
	$D(A)\subset X$. Let $Q\subset D(A)$. 	We define the {\it modulus of continuity }$\omega\left(\delta\right), \delta\in[0;\,\infty),$ of the operator $A$ on the class $Q$ 	to be
		\[
		\omega(\delta):=\sup\limits _{
			\begin{tiny}\begin{array}{c} x\in Q,\\
			\left\Vert x\right\Vert \leq\delta
			\end{array}\end{tiny}}
	\left\Vert Ax\right\Vert _{X}.
	\]
The problem of computing the modulus of continuity of an operator on the given class of elements is an abstract version of the problem of finding sharp Landau-Kolmogorov type inequality
(see, for instance,~\cite[Ch. 7]{babenko_inequalities}).

Let $\mathcal{\: L}(N)=\mathcal{L}(N;X,Y)\,$ be a set of linear bounded operators $T$ from $X$ to $Y$ with norms
	$\|T\|=\|T\|_{X\to Y}$ bounded by $N>0.$	The quantity
	\[
		U(T)=\sup\limits_{x\in Q}\|Ax-Tx\|_{Y}
	\]
	is called the {\it deviation} of the operator $T\in\mathcal{L}(N)$ from the operator $A$ on the class $Q$. Finally,
	\begin{equation}
		E(N)=E(N;A,Q):=\inf\limits_{T\in\mathcal{L}(N)}\{U(T)\label{eq:1}
	\end{equation}
	is called {\it the best approximation} of operator $A$ by a set of bounded operators $\mathcal{L}(N)$ on the class $Q.$
	
	Stechkin's problem (see, for instance, \cite{arestov_approx},~\cite{stechkin_inequalities},~\cite{stechkin_best_approx}, and~\cite{babenko_inequalities}, Ch. 7.1) of the best approximation of the operator $A$ on class $Q$ consists of computing $E(N)$ and finding (studying the questions of its existence, uniqueness, characterization) the ``extremal'' operator, i.e. the one that delivers $\inf$ in (\ref{eq:1}).

Let ${\cal O}={\cal O}(X,Y)$ be the set of all mappings of the space $X$ into the space $Y$,
 ${\cal L}={\cal L}(X,Y)$ be the set of all linear operators from $X$ into $Y$, and
${\cal B}={\cal B}(X,Y)$ be the set of all bounded linear operators $X$ into $Y$.

For $\delta\ge 0$ and operator $T\in {\cal R}$ we set
\[
U_\delta(T)=U_\delta(T;A,Q):=\sup\limits_{x\in Q, \eta\in X\atop \| x-\eta\|_X\le \delta}\| Ax-T\eta\|_Y;
\]
\[
{\cal E}_\delta({\cal R})={\cal E}_\delta({\cal R}; A,Q)=\inf\limits_{T\in {\cal R}}U_\delta(T).
\]
The problem of optimal recovery of the operator $A$ with the help of the set of mappings ${\cal R}$ (recovery methods) on the elements from the class $Q$, given with error $\delta$, consists of finding the quantity ${\cal E}_\delta({\cal R})$ and operator $T\in {\cal R}$, which realizes $\inf\limits_{T\in {\cal R}}U_\delta(T)$.

		Let $Q$ be a set in $X$ and $x\in X.$ The quantity
	\[
	E(x,Q)_{X}=\inf\limits _{q\in Q}\|x-q\|_{X},
	\]
	is called the {\it best approximation of the element $x\in X$} by the set $Q\subset X.$
	
	Let $F,Q$ be convex classes in $X$ and $N>0$ is a real number. The set $NQ=\{Nx:x\in Q\}$ is called a {\it homothet} of the class	$Q$ with the homothety coefficient $N.$ The quantity	
     \begin{equation}
	E(F,NQ)_{X}=\sup\limits _{u\in F}E(u,NQ)_{X}=\sup\limits _{u\in F}\inf\limits _{x\in NQ}\|u-x\|_{X}\label{eq:12}
	\end{equation}
	is called the {\it best approximation of the class $F$ by a homothet $NQ.$ }
	
	The problem of best approximation of the class $F$ by a homothet $NQ$ consist of computing the quantity (\ref{eq:12}).

 All above listed problems are closely related to Landau-Kolmogorov type inequalities (see, for instance,~\cite{arestov_approx},~\cite{babenko_inequalities} \S\S 7.3 - 7.5).
We need the following two theorems, which formally establish the connection.

Set
	\[
\Delta(N):=\sup\limits_{\delta >0}\{ \omega (\delta)-N\delta\}
\]
and
\[
l(\delta):=\inf\limits_{ N\ge 0}\{ E(N)+N\delta\}.
\]
	\begin{theorem}(Stechkin~\cite{stechkin_best_approx}).\label{Th1}
		Let $A$ be a homogeneous (in particular, linear) operator and let $Q$ be a centrally symmetric subset of $D(A)$. Then
		\[
			E(N)\geq \Delta(N), \qquad N\geq 0,
		\]
and		
		\[
			\omega(\delta)\leq l(\delta), \qquad \delta \geq 0.
		\]
	\end{theorem}

\begin{theorem}\label{Th2}
If $Q$ is a centrally symmetric set and $A$ is a homogeneous operator, then
\[
\omega(\delta)\le {\cal E}_\delta({\cal O})\le {\cal E}_\delta({\cal B})={\cal E}_\delta({\cal L})\le l(\delta).
\]
\end{theorem}

In this paper, we consider the above stated problems in the case when $X=Y=H$, where $H$ is a Hilbert space, the considered operators are some functions of self-adjoint operator in $H$, and the class of elements is also defined with the help of some function of the same operator.

Let us mention some known results of this type for operators in Hilbert spaces.

 First of all, let us mention the classical Hardy-Littlewood-Polya inequality~\cite{hardy_inequalities} for function $x(t)$ from $L_2(\mathbb R)$, such that the derivative $x^{(r)}(t)$ of order $r\in{\mathbb N}$ in Sobolev sense also belongs to the space $L_2(\mathbb R)$:
\begin{equation}\label{hlp}
\| x^{(k)}\|_{L_2({\mathbb R})}\le \| x\|_{L_2({\mathbb R})}^{1-\frac kr}\| x^{(r)}\|_{L_2({\mathbb R})}^{\frac kr},\qquad k=1,...,r-1.
\end{equation}

This inequality implies the following estimate for the modulus of continuity of the operator $\frac{d^k}{dt^k}$ on the class $W^r_{2,2}({\mathbb R})=\{ x\in L_2({\mathbb R})\; :\; \| x^{(r)}\|_{L_2({\mathbb R})}\le 1\}$:
\[
\omega(\delta)\le \delta^{1-\frac kr}.
\]
In addition, the sharpness of Hardy-Littlewood-Polya inequality implies that in fact
\[
\omega(\delta)= \delta^{1-\frac kr},\qquad \delta \ge 0.
\]

Hardy-Littlewood-Polya inequalty and the above result on computing of the modulus of continuity
have been generalized in multiple directions (see, for instance,~\cite{babenko_inequalities},~\cite{babenko_ligun_shumeiko},~\cite{babenko_approx},~\cite{Bab_Kryachko}). In particular, in~\cite{Bab_Kryachko} sharp Hardy-Littlewood-Polya type inequality was proved for functions of operators with a discrete spectrum.

The problem of best approximation of the unbounded operator $\frac{d^k}{dt^k}$ by bounded operators on the class $W^r_{2,2}({\mathbb R})$ was solved in
~\cite{subbotin_best_approx}. It was proved there that
\[
E(N)\le \frac{\frac{k}{r}\left(1-\frac{k}{r}\right)^{\frac{r-k}{k}}}{N^{\frac{r-k}{k}}}.
\]
The result was further generalized~\cite{babenko_approx} to the case of higher degrees of self-adjoint operators in Hilbert spaces.

In~\cite{subbotin_best_approx} it was also proved that
\[
E(W^{r-k}_{2,2}({\mathbb R}),NW^r_{2,2}({\mathbb R}))\le \frac{\frac{k}{r}\left(1-\frac{k}{r}\right)^{\frac{r-k}{k}}}{N^{\frac{r-k}{k}}}.
\]
This was generalized in~\cite{bapprclasses} to the case when classes are defined with the help of degrees of arbitrary self-adjoint operators.

The paper is organized as follows. In Section~\ref{S2} we introduce the necessary definitions and facts from spectral theory of self-adjoint operators in Hilbert spaces. In particular, here we define functions of such operators.
In Section~\ref{S3} we obtain rather general Hardy-Littlewood-Polya type inequality for functions of unbounded self-adjoint operators. In addition, we find the modulus of continuity of a function of an operator on a class of elements defined with the help of another function of an operator. In Section~\ref{S4} we solve the problem of best approximation of a function of an unbounded self-adjoint operator by bounded operators. In Section~\ref{S5} we obtain a series of sharp additive Hardy-Littlewood-Polya type inequalities for functions of operators. The problem of approximation of one class of elements from a Hilbert space by another class is solved in Section~\ref{S6}. Finally, in Section~\ref{S7} we solve the problem of optimal recovery of operators on a class of elements given with an error.

	\medskip{}
	
\section{Preliminaries from Spectral Theory}\label{S2}	

We begin by reminding some necessary facts on operator Stieltjes integrals and functions of self-adjoint operators in Hilbert spaces.


	Let $H$ be a Hilbert space with an inner product $\left(x,\, y\right)$ and norm
	$\left\Vert x\right\Vert =\left(x,x\right)^{1/2}.$ We consider a linear unbounded operator $A$ in $H$ with domain $D(A)$.

	First, let us recall some definitions and facts from spectral theory of self-adjoint operators (see, for instance, \S 75 and \S 88 in ~\cite{akhiezer_theory}).
	
{\it Partition of unity} is a one parametric family of projection operators $E_{t}$, defined on a finite or infinite interval $\left[\alpha,\,\beta\right]$ (in the case when the interval $\left[\alpha,\,\beta\right]$ is infinite, we understand, by definition,
	\[
		E_{-\infty}=\lim\limits _{t\rightarrow-\infty}E_{t},\qquad E_{\infty}=\lim\limits _{t\rightarrow\infty}E_{t}
	\]
	in strong convergence sense) and satisfying the following properties:
	
    a) $E_{u}E_{v}=E_{s},\; s=\min\left\{ u,\, v\right\} ,$
	
	b) in the sense of strong convergence
	\[
		E_{t-0}=E_{t}, \qquad \alpha<t<\beta,
	\]

	c) $E_{\alpha}=0,\; E_{\beta}=I$ ($I$ is an identity operator)
	
	We set $E_{t}=0$ for $t\leq\alpha$ and $E_{t}=I$ for $t\geq\beta.$
	
	It follows from the definition that for any $x\in H$ the quantity
	\[
		\sigma\left(t\right)=\left(E_{t}x,\, x\right),\qquad \alpha<t<\beta,
	\]
	is left-continuous, non-decreasing function of bounded variation for which
	\[
		\sigma\left(\alpha\right)=0,\qquad \sigma\left(\beta\right)=\left(x,\, x\right).
	\]
	Thus, we have $\sigma$-measure that allows the construction of Lebesgue-Stieltjes integral.
	
	If any condition is satisfied with respect to all $\sigma$-measures, generated by elements $x\in H$, then we say that it is satisfied with respect to the operator measure $E_t$.
	
	Now for the defined, measurable, and finite almost everywhere with respect to the operator measure $E_t$ functions, we may consider operator integrals (for details and properties of such integrals see, for instance, ~\cite{akhiezer_theory})
	$$
	\displaystyle \int_{-\infty}^{\infty} \varphi(t)dE_t.
	$$
	
	Based on the spectral theorem, each self-adjoint operator $A$ has a corresponding partition of unity $E_{t},\: t\in\mathbb{R}$,
	such that
	\[
		A=\int\limits _{-\infty}^{+\infty}tdE_{t}.
	\]
	In addition, element $x$ belongs to the domain $D\left(A\right)$ of the operator $A$ if and only if
	\[
		\int\limits _{-\infty}^{+\infty}t^{2}d\left(E_{t}x,x\right)<\infty.
	\]
	Moreover, if $x\in D\left(A\right),$ then
	\[
		Ax=\int\limits _{-\infty}^{+\infty}tdE_{t}x,
	\]
	and
	\[
		\left\Vert Ax\right\Vert ^{2}=\int\limits _{-\infty}^{+\infty}t^{2}d\left(E_{t}x,x\right)<\infty.
	\]
Let now function $\varphi(t)$ be defined, measurable, and finite almost everywhere with respect to the operator measure $E_t$. We also assume that there exists a dense set $D$ in $H$ of elements $x$, such that
\begin{equation}\label{dom}
\int\limits _{-\infty}^{+\infty}|\varphi(t)|^{2}d\left(E_{t}x,x\right)<\infty.
\end{equation}
Under the made assumptions, the function $\varphi(A)$ of an operator $A$ is an operator defined as follows
\[
\varphi(A)x=\int\limits _{-\infty}^{+\infty}\varphi(t)dE_{t}x
\]
for all those $x\in H$ such that (\ref{dom}) holds. Relation (\ref{dom}) defines the domain $D(\varphi(A))$ of the operator $\varphi(A)$.

	\section{Inequalities of Hardy-Littlewood-Polya type and the problem of computing the modulus of continuity}\label{S3}
	We begin with the case when rather general Hardy-Littlewood-Polya type inequality can be proved in a simple and explicit manner.
	We consider functions $\varphi(A)$ and $\psi(A)$ of an operator $A$, where $\varphi(t)$ and $\psi\left(t\right)$ are continuous complex-valued functions on $\mathbb{R},$ such that $|\varphi(t)|$ and $|\psi(t)|$ are even and strictly increasing on $(0,\infty)$. In addition, we assume
	\begin{equation}\label{connect}
		|\varphi(t)|^{2}=F(|\psi(t)|^{2})
	\end{equation}
	where $F(\cdot)$ is a strictly increasing, concave function, and $F(0)=0.$
	
A rather general  Hardy-Littlewood-Polya type inequality is contained in the following theorem.

	\begin{theorem}\label{Th3}
		Let $A$ be an unbounded self-adjoint operator in $H$. Then for any $x\in D(\psi(A)),\; x\neq\theta ,$ the following inequality holds
		\begin{equation}
			\left\Vert \varphi(A)x\right\Vert^2 \leq\left\Vert x\right\Vert^2 {F\left(\frac{\left\Vert \psi(A)x\right\Vert ^{2}}{\left\Vert x\right\Vert ^{2}}\right)}.\label{eq:4}
		\end{equation}
If, in addition, $A$ is such that
		\begin{equation}
			\left(E_{t}-E_{s}\right)D\left(\psi(A)\right)\neq\left\{ \theta\right\} ,\qquad 0\leq s<t\leq\infty,\label{eq:3}
		\end{equation}
		then the inequality (\ref{eq:4}) is exact.
	\end{theorem}

{\bf Remark.} КThe classical Hardy-Littlewood-Polya inequality in the multiplicative form (\ref{hlp}) can be obtained by taking $A=i\frac{d}{dt}$, $\varphi(t)=t^k$, $\psi(t)=t^r$, and $F(t)=t^{\frac kr}$. Hence, we call the form in (\ref{eq:4}) the {\it multiplicative form}.

\proof
 In order to obtain the upper estimate, we apply Jensen's inequality ($F$ is a concave function and $\int\limits _{-\infty}^{+\infty}\frac{d\left(E_{t}x,x\right)}{\left\Vert x\right\Vert ^{2}}=1$).
		For any $x\in D\left(\psi(A)\right),\; x\neq\theta, $ we obtain
		\[
		\begin{array}{lll}
			\left\Vert \varphi(A)x\right\Vert ^{2}&=&\displaystyle \int\limits _{-\infty}^{+\infty}|\varphi(t)|^2{d\left(E_{t}x,x\right)}=\left\Vert x\right\Vert ^{2}\displaystyle \int\limits _{-\infty}^{+\infty}F\left(\psi^{2}\left(t\right)\right)\frac{d\left(E_{t}x,x\right)}{\left\Vert x\right\Vert ^{2}}\cr
			&\leq& \displaystyle \left\Vert x\right\Vert ^{2}F\left(\frac{1}{\left\Vert x\right\Vert ^{2}}\displaystyle\int\limits _{-\infty}^{+\infty}\left|\psi\left(t\right)\right|^{2}d\left(E_{t}x,x\right)\right)\cr
			&=&\displaystyle \left\Vert x\right\Vert ^{2}F\left(\frac{\left\Vert \psi(A)x\right\Vert ^{2}}{\left\Vert x\right\Vert ^{2}}\right).
			\end{array}
		\]
		Therefore,
		\begin{equation}\nonumber
			\left\Vert \varphi(A)x\right\Vert^2 \leq\left\Vert x\right\Vert^2 {F\left(\frac{\left\Vert \psi(A)x\right\Vert ^{2}}{\left\Vert x\right\Vert ^{2}}\right)}.
		\end{equation}

The fact that under assumption (\ref{eq:3}) obtained inequality is sharp will follow from the next theorem.
$\square$
	
	By $W^{\psi}$ we denote the class of elements $x\in D\left(\psi\left(A\right)\right)$ such that $\left\Vert \psi(A)x\right\Vert \leq1.$
		
	\begin{theorem}\label{Th4}
		Let $A$ be an unbounded self-adjoint operator in $H$ and let $\omega(\delta)$ be the modulus of continuity of the operator
		$\varphi(A)$ on the class $W^{\psi}$.
		Then for any $\delta>0$
		\[
			\omega(\delta)\leq\delta\sqrt{F\left(\frac{1}{\delta^{2}}\right)}.
		\]
		If, in addition, $A$ is such that the assumption (\ref{eq:3}) is satisfied,
				then for any $\delta>0$
		\[
			\omega(\delta)=\delta\sqrt{F\left(\frac{1}{\delta^{2}}\right)}.
		\]
	\end{theorem}
	
	\proof
				Inequality~(\ref{eq:4}), together with assumptions imposed on the function $F$ and the fact that $x \in W^{\psi}$ and $\|x\|\leq \delta$, implies
		\[
			\left\Vert \varphi(A)x\right\Vert^2 \leq\delta F\left(\frac{1}{\delta^{2}}\right),
		\]
		and, hence,
				\begin{equation}
			\omega(\delta)\leq\delta\sqrt{F\left(\frac{1}{\delta^{2}}\right)}.\label{eq:5}
		\end{equation}

		Next, let us show that if for the operator $A$ assumption (\ref{eq:3}) is satisfied, then we also have the following lower estimate
		\[
			\omega(\delta)\geq\delta\sqrt{F\left(\frac{1}{\delta^{2}}\right)}.
		\]
		For a given $\delta>0$ and arbitrary $\varepsilon\in\left(0,\,1\right)$, we take
		\[
			t=|\psi|^{-1}\left(\frac{1}{\delta}\right)\;\;\;\;\text{and}\;\;\;\; s=(1-\varepsilon)|\psi|^{-1}\left(\frac{1}{\delta}\right),
		\]
		where $|\psi|^{-1}(t)$ is an inverse function to $|\psi|(t).$ We choose an element $x_{\delta,\varepsilon}\in\left(E_{t}-E_{s}\right)D\left(\psi(A)\right)$ in such a way that $\left\Vert x_{\delta,\varepsilon}\right\Vert =\delta.$
		
		First of all, it is easy to check that the chosen element $x_{\delta,\varepsilon}$ belongs to the class $W^{\psi}.$
		Indeed, since $|\psi(u)|$ is strictly increasing, we have
		\[
		\begin{array}{lll}
			\left\Vert \psi(A)x_{\delta,\varepsilon}\right\Vert ^{2}&=&\displaystyle \int\limits _{s}^{t}\left|\psi(u)\right|^{2}d\left(E_{u}x_{\delta,\varepsilon},x_{\delta,\varepsilon}\right)\cr
			&\leq&
		|\psi|^{2}(t)\left\Vert x_{\delta,\varepsilon}\right\Vert ^{2}=\left(|\psi|(|\psi|^{-1}(\frac{1}{\delta}))\right)^{2}\delta^{2}=1.
		\end{array}
		\]

		Next, we apply operator $\varphi\left(A\right)$ to the element $x_{\delta,\varepsilon}$ and, taking into account the fact that $|\varphi(u)|$ is strictly increasing we obtain the following estimate from below
		\[
		\begin{array}{lll}
			\omega(\delta)^2\ge\left\Vert \varphi(A)x_{\delta,\varepsilon}\right\Vert ^{2}&=& \displaystyle \int\limits _{s}^{t}|\varphi|^{2}(u)d\left(E_{u}x_{\delta,\varepsilon},x_{\delta,\varepsilon}\right)\geq\varphi^{2}(s)\left\Vert x_{\delta,\varepsilon}\right\Vert ^{2}\cr
			&=&|\varphi|^{2}((1-\varepsilon)t)\delta^{2}=F(|\psi|^{2}\left((1-\varepsilon)|\psi|^{-1}(\frac{1}{\delta})\right))\delta^{2}.
			\end{array}
		\]
		As $\varepsilon\rightarrow0$, we arrive at
		\begin{equation}
			\omega(\delta)\ge\sqrt{F(|\psi|^{2}\left(|\psi|^{-1}(\frac{1}{\delta})\right))\delta^{2}} =\delta\sqrt{F\left(\frac{1}{\delta^{2}}\right)}.\label{eq:6}
		\end{equation}
		Combining (\ref{eq:5}) and (\ref{eq:6}), we obtain
		\[
			\omega(\delta)=\delta\sqrt{F\left(\frac{1}{\delta^{2}}\right)}.
		\]
		$\square$

	\section{The problem of best approximation of an unbounded operator by bounded ones}\label{S4}
In this section we consider the problem of approximating the operator $\varphi(A)$ on the class $W^\psi$ by bounded operators.

Let continuous, complex valued functions $\varphi(t)$ and $\psi(t)$ be such that $|\varphi(t)|$ and $|\psi(t)|$ are even and strictly increasing to $+\infty$ on $(0,+\infty)$.
	For any $b>0$ let function $\varphi_{b}(t)$ be defined as follows
		\begin{equation}\label{extr}
			\varphi_{b}(t)=\begin{cases}
			\varphi(t)-\frac{\varphi(b)}{\psi(b)}\psi(t), & \left|t\right|\leq b\\
			0, & \left|t\right|\geq b.
			\end{cases}
		\end{equation}
Let also
		\[
			\max_{t}\left|\varphi_{b}(t)\right|=N(b).
		\]

	\begin{theorem}\label{Th5}
		Let $A$ be an unbounded, self-adjoint operator in a Hilbert space $H.$ Let also functions $\varphi$ and $\psi$ be such that the function $\frac{|\phi(t)|}{|\psi(t)|}$ is non-increasing. Then for any $b>0$
			\begin{equation}
			E\left(N(b)\right)\leq\frac{|\varphi(b)|}{|\psi(b)|}.\label{eq:7}
		\end{equation}
				If functions $\varphi$ and $\psi$ are such that (\ref{connect}) holds and operator $A$ is such that condition (\ref{eq:3}) is satisfied, then for any $b>0$
			\begin{equation}
			E\left(N(b)\right)=\frac{|\varphi(b)|}{|\psi(b)|}= \frac 1{|\psi(b)|}\sqrt{F(|\psi(b)|^2)}=\omega\left(\frac 1{|\psi(b)|}\right).\label{eq:8}
		\end{equation}
where $\omega(\delta)$ is modulus of continuity of the operator $\varphi$ on the class $W^\psi$.

				The extremal operator is $\varphi_{b}(A)$ defined with the help of a function in (\ref{extr}).

If, additionally, the function $\frac{|\phi(t)|}{|\psi(t)|}$ on the interval $(0, +\infty)$ is strictly decreasing from $+\infty$ to $0$, then for any $N>0$
\[
E(N)=\frac{|\varphi(b)|}{|\psi(b)|},
\]
where $b$ is the unique solution of the equation $N(b)=N.$

 	\end{theorem}

\proof
		With the help of functions $\varphi_{b}(t)$, we define operators $\varphi_{b}(A)$ and estimate their norms $\left\Vert \varphi_{b}(A)\right\Vert :$
		\[
			\left\Vert \varphi_{b}(A)x\right\Vert ^{2}=\int\limits _{-b}^{b}|\varphi(u)|^{2}d(E_{u}x,x)\leq\max_{t}\left|\varphi_{b}(t)\right|^{2}\left\Vert x\right\Vert ^{2},
		\]
		from which it follows that
		\[
			\left\Vert \varphi_{b}(A)\right\Vert \leq\max_{t}\left|\varphi_{b}(t)\right| = N(b).
		\]

		Next, we find the upper estimate for the norm $\left\Vert \varphi(A)x-\varphi_{b\left(N\right)}(A)x\right\Vert.$ For $x\in W^{\psi}$ we have
		\[
		\begin{array}{lll}
			\left\Vert \varphi(A)x-\varphi_{b}(A)x\right\Vert ^{2}&=&\displaystyle \int\limits _{-\infty}^{+\infty}\left|\varphi(u)-\varphi_{b}(u)\right|{}^{2}d(E_{u}x,x)\cr
					&\leq& \displaystyle \max_{t}\frac{\left|\varphi(t)-\varphi_{b}(t)\right|{}^{2}}{|\psi(t)|^{2}}\int\limits _{-\infty}^{+\infty}\left|\psi(u)\right|^{2}d(E_{u}x,x)\cr
					&=&\displaystyle \max_{t}\frac{\left|\varphi(t)-\varphi_{b}(t)\right|{}^{2}}{|\psi(t)|^{2}}\left\Vert \psi(A)x\right\Vert ^{2} \leq \displaystyle \max_{t}\frac{\left|\varphi(t)-\varphi_{b}(t)\right|{}^{2}}{|\psi(t)|^{2}}.
					\end{array}
		\]
		Taking into account the fact that the function $\frac{|\varphi\left(y\right)|^{2}}{|\psi\left(y\right)|^{2}}$ is non-increasing, we obtain		
		\[
			\frac{|\varphi(t)-\varphi_{b}(t)|^{2}}{|\psi(t)|^{2}}=\begin{cases}
			\frac{|\varphi(b)|^{2}}{|\psi(b)|^{2}}, & \left|t\right|\leq b\\
			\frac{|\varphi(t)|^{2}}{|\psi(t)|^{2}}, & \left|t\right|\geq b
			\end{cases}
		\]
		 and
		\[
			\max_{t}\left|\frac{\varphi(t)-\varphi_{b}(t)}{\psi(t)}\right|=\frac{|\varphi(b)|}{|\psi(b)|}.
		\]
		From here it follows that for $x\in W^{\psi}$
		\begin{equation}
			E\left(N(b)\right)\leq\sup\limits _{
				\begin{tiny}\begin{array}{c}x\in W^{\psi}\end{array}\end{tiny}}
			\left\Vert \varphi(A)x-\varphi_{b}(A)x\right\Vert \leq\frac{|\varphi(b)|}{|\psi(b)|}.\label{eq:10}
		\end{equation}
		Inequality (\ref{eq:7}) is proved.
		
		If, in addition, condition (\ref{connect}) and assumption (\ref{eq:3}) are satisfied, then by Theorem~\ref{Th4}
		\[
			\omega(\delta)=\delta\sqrt{F\left(\frac{1}{\delta^{2}}\right)}.
		\]
		Taking into account Theorem~\ref{Th1}, we obtain
		\[
		\begin{array}{lll}
			E\left(N(b)\right)&\geq&\displaystyle \sup_{\delta\geq0}\left\{ \omega\left(\delta\right)-N(b)\delta\right\} \cr
					&=&\displaystyle \sup_{\delta\geq0}\left\{ \delta\sqrt{F\left(\frac{1}{\delta^{2}}\right)}-\delta\max\limits_t\left|\varphi_{b}(t)\right|\right\}.
\end{array}
\]
Let $\xi$ be such that
\[
\max\limits_t|\varphi_b(t)|=|\varphi_b(\xi)|.
\]
Then
\[
E\left(N(b)\right)\ge \sup_{\delta\geq0}\left\{ \delta\sqrt{F\left(\frac{1}{\delta^{2}}\right)}-\delta\left|\varphi(\xi)-
\frac{\varphi(b)}{\psi(b)}\psi(\xi)\right|\right\}.
\]
Taking $\delta=\frac{1}{|\psi\left(\xi\right)|},$ we find
		\begin{equation}
			E\left(N(b)\right)\geq \frac{|\varphi(\xi)|}{|\psi(\xi)|}-\left|\frac{\varphi(\xi)}{\psi(\xi)}- \frac{\varphi(b)}{\psi(b)}\right| \geq \frac{|\varphi(b)|}{|\psi(b)|}. \label{eq:11}
		\end{equation}
		Combining (\ref{eq:10}) and (\ref{eq:11}), we obtain the desired (\ref{eq:8}).

Finally, let us prove the last statement of the theorem.
 We prove that under the made assumptions, for any $N>0$ there exists $b>0$ such that $N(b)=N$.
Clearly, $N(b)$ continuously depends on $b$ and is strictly increasing with $b$. Therefore, it is sufficient to prove that for any $M>0$ there exists $b>0$ such that $N(b)>b$.

Let us take an arbitrary $M>0$. Since $|\varphi(t)|\to +\infty$ as $t\to +\infty$ we then choose $c>0$ so that $|\varphi(c)|>2M$. In addition,
\begin{equation}\label{ots}
N(b)=\max\limits_t|\varphi_b(t)|\ge |\varphi_b(c)|\ge |\varphi (c)|-\frac{|\varphi (b)|}{|\psi (b)|}|\psi(c)|.
\end{equation}
Since $\frac{\varphi(t)}{\psi(t)}\to 0 $ as $t\to +\infty$, for all large enough $b$ we have
\[
\frac{|\varphi (b)|}{|\psi (b)|}|\psi(c)|<M.
\]
Taking into account (\ref{ots}), we obtain
\[
N(b)\ge 2M -M=M,
\]
as desired. 
$\square$

\section{Some additive inequalities of Hardy-Littlewood-Polya type }\label{S5}

Recall that for any $b>0$ the function $\varphi_{b}(t)$ is defined by (\ref{extr}). As before, let
		\[
			\max_{t}\left|\varphi_{b}(t)\right|=N(b).
		\]

	\begin{theorem}\label{Th6}
		Let $A$ be a self-adjoint operator in a Hilbert space $H$ and let functions $\varphi$ and $\psi$ be such that the function $\frac{|\varphi(t)|}{|\psi(t)|}$ is non-increasing. Then for any $b>0$ and any
		$x\in D\left(\psi\left(A\right)\right)$ we have
		\begin{equation}
		\left\Vert \varphi(A)x\right\Vert \leq\frac{|\varphi(b)|}{|\psi(b)|}\left\Vert \psi(A)x\right\Vert +N(b)\left\Vert x\right\Vert .\label{eq:15}
		\end{equation}
		If the operator $A$ is such that condition (\ref{eq:3}) is satisfied,
		then for any $b>0$ the inequality (\ref{eq:15}) is sharp in the sense that the constant $\frac{|\varphi(b)|}{|\psi(b)|}$ cannot be taken smaller.\end{theorem}
		
	\proof
		For $x\in D\left(\psi\left(A\right)\right)$, using the triangle inequality, we have
		
		\begin{equation}
		\left\Vert \varphi(A)x\right\Vert \leq\left\Vert \varphi(A)x-\varphi_{b}(A)x\right\Vert +\left\Vert \varphi_{b}(A)x\right\Vert .\label{eq:17}
		\end{equation}
		While proving Theorem~\ref{Th5} we obtained the following estimate from above for the norm of the difference $\left\Vert \varphi(A)x-\varphi_{b}(A)x\right\Vert$:
		\[
		\left\Vert \varphi(A)x-\varphi_{b}(A)x\right\Vert \leq\max_{t}\left|\frac{\varphi(t)-\varphi_{b}(t)}{\psi(t)}\right|\left\Vert \psi(A)x\right\Vert ,
		\]
		or, taking into account the equality,
		\[
		\max_{t}\left|\frac{\varphi(t)-\varphi_{b}(t)}{\psi(t)}\right|=\frac{|\varphi(b)|}{|\psi(b)|},
		\]
		becomes
		\begin{equation}
		\left\Vert \varphi(A)x-\varphi_{b}(A)x\right\Vert \leq\frac{|\varphi(b)|}{|\psi(b)|}\left\Vert \psi(A)x\right\Vert .\label{eq:18}
		\end{equation}
		Next, we consider $\left\Vert \varphi_{b}(A)x\right\Vert .$
		Before, we obtained the inequality
		\[
		\left\Vert \varphi_{b}(A)x\right\Vert \le \max\limits_{t}\left|\varphi_{b}(t)\right|\left\Vert x\right\Vert =N(b)\left\Vert x\right\Vert ,
		\]
and, therefore,
		\begin{equation}
		\left\Vert \varphi_{b}(A)x\right\Vert \leq N(b)\left\Vert x\right\Vert .\label{eq:19}
		\end{equation}
	Combining (\ref{eq:18}), (\ref{eq:19}), and (\ref{eq:17}), we obtain (\ref{eq:15}).
		
		Next, we show that if assumption (\ref{eq:3}) is satisfied, then the inequality
		(\ref{eq:15}) is sharp.
		Assume to the contrary, that there exists $\delta>0$ such that for any $x\in D(\psi(A))$
		\[
		\left\Vert \varphi(A)x\right\Vert \leq\left(1-\delta\right)\frac{|\varphi(b)|}{|\psi(b)|}\left\Vert \psi(A)x\right\Vert +N(b)\left\Vert x\right\Vert .
		\]
		Let $\xi$ be such that
\[
\max\limits_t|\varphi_b(t)|=|\varphi_b(\xi)|.
\]
For arbitrary
		$\varepsilon\in\left(0,\,1\right)$, we set $t=\xi$
		and $s=\left(1-\varepsilon\right)t$. Let $x_{\xi,\varepsilon}\in\left(E_{t}-E_{s}\right)D\left(\psi(A)\right)$ and $x_{\xi,\varepsilon}\neq \theta$. Then
		\[
		\left\Vert \varphi(A)x_{\xi,\varepsilon}\right\Vert ^{2}=\int\limits _{s}^{t}\left|\varphi(u)\right|^{2}d(E_{u}x_{\xi,\varepsilon},x_{\xi,\varepsilon})\geq|\varphi|^{2}\left(s\right)\left\Vert x_{\xi,\varepsilon}\right\Vert ^{2}=|\varphi|^{2}\left(\left(1-\varepsilon\right)t\right)\left\Vert x_{\xi,\varepsilon}\right\Vert ^{2}
		\]
		and
		\[
		\left\Vert \psi(A)x_{\xi,\varepsilon}\right\Vert ^{2}=\int\limits _{s}^{t}\left|\psi(u)\right|^{2}d(E_{u}x_{\xi,\varepsilon},x_{\xi,\varepsilon})\leq |\psi|^{2}\left(t\right)\left\Vert x_{\xi,\varepsilon}\right\Vert ^{2}.
		\]	
Taking into account these inequalities, we have		
		$$
		\begin{array}{lll}
		|\varphi|\left(\left(1-\varepsilon\right)t\right)\left\Vert x_{\xi,\varepsilon}\right\Vert &\leq& \displaystyle \left\Vert \varphi(A)x_{\xi,\varepsilon}\right\Vert \leq\left(1-\delta\right)\frac{|\varphi(b)|}{|\psi(b)|}\left\Vert \psi(A)x_{\xi,\varepsilon}\right\Vert +N(b)\left\Vert x_{\xi,\varepsilon}\right\Vert \cr
		&\leq& \displaystyle \left(1-\delta\right)\frac{|\varphi(b)|}{|\psi(b)|}\psi\left(t\right)\left\Vert x_{\xi,\varepsilon}\right\Vert +N(b)\left\Vert x_{\xi,\varepsilon}\right\Vert .
		\end{array}
		$$
Therefore,
		\[
		|\varphi|\left(\left(1-\varepsilon\right)t\right)\leq\left(1-\delta\right)\frac{|\varphi(b)|}{|\psi(b)|}\psi\left(t\right)+N(b).
		\]	
		Since $\varepsilon$ is arbitrary, we have
		\[
		|\varphi\left(t\right)|\leq\left(1-\delta\right)\frac{|\varphi(b)|}{|\psi(b)|}|\psi|\left(t\right)+N(b)
		\]
		or
		\[
		\frac{|\varphi|\left(t\right)-N(b)}{|\psi|\left(t\right)}\leq\left(1-\delta\right)\frac{|\varphi(b)|}{|\psi(b)|},
		\]
		which, together with $\left|\varphi_{b}(\xi)\right|=\left|\varphi(\xi)-\frac{\varphi(b)}{\psi(b)}\psi(\xi)\right|=N(b)$
		for the chosen $t=\xi,$ implies that
		\[
		\frac{|\varphi|\left(\xi\right)-\left| \varphi(\xi)-\frac{\varphi(b)}{\psi(b)}\psi(\xi)\right| }{|\psi\left(\xi\right)|}\leq\left(1-\delta\right)\frac{|\varphi(b)|}{|\psi(b)|},
		\]
		i.e.
		
		\[
		1\leq\left(1-\delta\right).
		\]
		Since $\delta>0,$ the last inequality is not possible and we obtained the desired contradiction. Therefore, inequality (\ref{eq:15}) is sharp. $\square$

\begin{cor}
	If the operator $A$ is such that condition (\ref{eq:3}) is satisfied
		then for any $b>0$
\[
	\sup_{x\in D(\psi(A))}\frac{\left\Vert \varphi(A)x\right\Vert -N(b)\left\Vert x\right\Vert}{\left\Vert \varphi(A)x\right\Vert}	 =\frac{|\varphi(b)|}{|\psi(b)|}.
\]
    \end{cor}

	\section{Best approximation of some class of elements of a Hilbert space by another class.}\label{S6}
	
	The next theorem provides the solution to the problem of best approximation of the class
	$F=W^{\frac{\psi}{\varphi}}$ by the homothet $NQ=NW^{\psi}.$
	\begin{theorem}\label{Th7}
		Let $A$ be a self-adjoint, unbounded operator in a Hilbert space $H.$ Then for any $b>0$
		
		\begin{equation}
		E(W^{\frac{\psi}{\varphi}},N(b)W^{\psi})\leq\frac{|\varphi(b)|}{|\psi(b)|}.\label{eq:20}
		\end{equation}
If the operator $A$ is such that condition (\ref{eq:3}) is satisfied,
		then for any $b>0$
		\begin{equation}
		E(W^{\frac{\psi}{\varphi}},N(b)W^{\psi})=\frac{|\varphi(b)|}{|\psi(b)|}.\label{eq:21}
		\end{equation}
	\end{theorem}
	\proof
Let
		\[
		\eta_{b}\left(t\right)=\begin{cases}
		1-\frac{\varphi(b)}{\psi(b)}\frac{\psi(t)}{\varphi\left(t\right)}, & \left|t\right|\leq b\\
		0, & \left|t\right|\geq b
		\end{cases}.
		\]
		Let us take an arbitrary element $x\in W^{\frac{\psi}{\varphi}}$ and consider $\eta_{b}\left(A\right)x.$
		We have
		$$
		\begin{array}{lll}
		\left\Vert \psi\left(A\right)\eta_{b}\left(A\right)x\right\Vert ^{2}&=&\displaystyle \int\limits _{-\infty}^{+\infty}\left|\psi\left(t\right)\right|^{2}\left|\eta_{b}(t)\right|^{2}d(E_{t}x,x)\cr
		&=&\displaystyle \int\limits _{-\infty}^{+\infty}\left|\varphi\left(t\right)\right|^{2}\left|\eta_{b}(t)\right|^{2}\cdot\frac{\left|\psi\left(t\right)\right|^{2}}{\left|\varphi\left(t\right)\right|^{2}}d(E_{t}x,x)\cr
		&\leq&\displaystyle \max_{t}\left|\varphi\left(t\right)\eta_{b}(t)\right|^{2}\left\Vert \frac{\psi}{\varphi}\left(A\right)x\right\Vert ^{2}
		\leq\displaystyle \max_{t}\left|\varphi_{b}(t)\right|^{2}=N(b)^{2}.
		\end{array}
		$$
		Thus,
		\[
		\left\Vert \psi\left(A\right)\eta_{b}\left(A\right)x\right\Vert \leq N(b),
		\]
		which implies $\eta_{b}\left(A\right)x\in N(b)W^{\psi}.$
		
		Next, for $x\in W^{\frac{\psi}{\varphi}}$ we find the estimate from above for the norm of the difference $\left\Vert x-\eta_{b}\left(A\right)x\right\Vert$:
		$$
		\begin{array}{lll}
		\left\Vert x-\eta_{b}\left(A\right)x\right\Vert ^{2}&=&\displaystyle \int\limits _{-\infty}^{+\infty}\left|1-\eta_{b}(t)\right|^{2}d\left(E_{t}x,x\right)\cr
				&\leq&\displaystyle \int\limits _{-\infty}^{+\infty}\frac{|\varphi|^{2}(b)}{|\psi|^{2}(b)}\frac{|\psi|^{2}(t)}{|\varphi|^{2}\left(t\right)}d\left(E_{t}x,x\right)=\frac{|\varphi|^{2}(b)}{|\psi|^{2}(b)}\left\Vert \frac{\psi}{\varphi}\left(A\right)x\right\Vert ^{2}\cr
				&\leq& \displaystyle \frac{|\varphi|^{2}(b)}{|\psi|^{2}(b)}.
		\end{array}
		$$
Therefore,
		\begin{equation}
		E(W^{\frac{\psi}{\varphi}},N(b)W^{\psi})\leq\frac{|\varphi(b)|}{|\psi(b)|}.\label{eq:22}
		\end{equation}
		On the other hand,
$$
\begin{array}{lll}
		E(W^{\frac{\psi}{\varphi}},N(b)W^{\psi})&=&\displaystyle \sup_{x\in W^{\psi/\varphi}}\inf_{y\in N(b)W^{\psi}}\sup_{\| u\|\le 1}\left|\left(u,\, x\right)-\left(u,\, y\right)\right|\cr
		&\geq & \displaystyle \sup_{x\in W^{\psi/\varphi}}\sup_{v\in W^\psi}\inf_{y\in N(b)W^{\psi}}\left|\left(\psi\left(A\right)v,\, x\right)-\left(\psi\left(A\right)v,\, y\right)\right|\cr
&\geq & \displaystyle \sup_{x\in W^{\psi/\varphi}}\sup_{v\in W^\psi}\inf_{y\in N(b)W^{\psi}}\left[\left|\left(\psi\left(A\right)v,\, x\right)\right|-\left|\left(\psi\left(A\right)v,\, y\right)\right|\right]\cr
		&=&\displaystyle \sup_{v\in W^\psi}\left[\sup_{x\in W^{\psi/\varphi}}\left|\left(\varphi\left(A\right)v,\,\frac{\psi}{\varphi}\left(A\right)x\right)\right|-N(b)\sup_{y\in W^\psi}\left|\left(v,\,\psi\left(A\right)y\right)\right|\right].
		\end{array}
		$$
Since the sets $\left\{\frac\psi\varphi(A)x\; :\; x\in W^{\frac\psi\varphi}\right\}$ and $\left\{\psi(A)x\; :\; x\in W^\psi\right\}$ are dense in the unit ball of $H$, we obtain
\[
E(W^{\frac{\psi}{\varphi}},N(b)W^{\psi})\ge
\displaystyle \sup_{\left\Vert \psi\left(A\right)v\right\Vert \leq1}\displaystyle \left(\left\Vert \varphi\left(A\right)v\right\Vert -N(b)\left\Vert v\right\Vert \right).
\]

		If condition (\ref{eq:3}) is satisfied for the operator $A$, then by Corollary 1 we obtain
		\begin{equation}
		E(W^{\frac{\psi}{\varphi}},N(b)W^{\psi})\ge\frac{|\varphi(b)|}{|\psi(b)|}.\label{eq:23}
		\end{equation}
		From (\ref{eq:22}) and (\ref{eq:23}), the equality (\ref{eq:21}) follows, which completes the proof of the theorem. $\square$

\section{Optimal recovery of operators}\label{S7}

In this section we present the solution to the problem of optimal recovery of an operator in the case when $X=Y=H$, the operator is $\varphi(A)$ (where $A$ is a self-adjoint operator in $H$), on the class $W^\psi$ and ${\cal R}={\cal O}, {\cal B}, {\cal L}$. More precisely, we prove the following theorem.
\begin{theorem}\label{Th8} Let error $\delta\geq 0$ be given.
If ${\cal R}={\cal O}, {\cal B}, {\cal L}$ and operator $A$ is such that (\ref{eq:3}) is satisfied, then
\[
{\cal E}_\delta({\cal R};\varphi(A),W^\psi)=\omega(\delta)=\delta\sqrt{F\left(\frac 1{\delta^2}\right)},
\]
where $\omega(\delta)$ is the modulus of continuity of the operator $\varphi(A)$ on the class $W^\psi$.
\end{theorem}

\proof
 Due to Theorem~\ref{Th2}, in order to prove Theorem~\ref{Th8} it is sufficient to prove that
\[
\omega(\delta)=l(\delta).
\]

From Theorems~\ref{Th1} and~\ref{Th5} it follows
\[
\omega(\delta)\le \inf\limits_{b>0}\left\{\frac{|\varphi(b)|}{|\psi(b)|}+N(b)\delta\right\}.
\]
Let us show that if (\ref{eq:3}) is satisfied, then the last inequality holds with equality sign.

We choose an arbitrary $\delta>0$ and $\varepsilon\in (0,1)$. Set $t=|\psi|^{-1}\left(\frac 1\delta\right)$ and $s=(1-\varepsilon)t.$ Then we choose an element $x_{\delta,\varepsilon}\in (E_t-E_s)D(\psi(A))$ such that $\|x_{\delta,\varepsilon}\|=\delta$. As it was shown in the process of proving Theorem~\ref{Th4}, $\|\psi(A)x_{\delta,\varepsilon}\|\le1$ and $\|\varphi(A)x_{\delta,\varepsilon}\|\ge |\varphi((1-\varepsilon)t)|\delta.$ 
Therefore
\[
\omega(\delta)\ge \|\varphi(A)x_{\delta,\varepsilon}\|\ge |\varphi((1-\varepsilon)t)|\delta.
\]
Letting $\varepsilon$ tend to zero we obtain that for any $b>0$
\[
\begin{array}{lll}
\omega(\delta)\ge |\varphi(t)|\delta
&=&\displaystyle \frac{|\varphi(b)|}{|\psi(b)|}+\delta\left( |\varphi(t)|-\frac{|\varphi(b)|}{|\psi(b)|}\frac 1\delta\right)\cr
&=&\displaystyle \frac{|\varphi(b)|}{|\psi(b)|}+\delta\left( |\varphi(t)|-\frac{|\varphi(b)|}{|\psi(b)|}|\psi(t)|\right).
\end{array}
\]
With an appropriate choice of $b=b_t$ we have
\[
|\varphi(t)|-\frac{|\varphi(b)|}{|\psi(b)|}|\psi(t)|=N(b_t),
\]
so that
\[
\begin{array}{lll}
\omega(\delta)& \ge & \displaystyle \frac{|\varphi(b_t)|}{|\psi(b_t)|}+\delta\left( |\varphi(t)|-\frac{|\varphi(b_t)|}{|\psi(b_t)|}|\psi(t)|\right)\cr
&=&\displaystyle \frac{|\varphi(b_t)|}{|\psi(b_t)|}+\delta N(b_t)\ge \inf\limits_{b>0}\left\{\frac{|\varphi(b)|}{|\psi(b)|}+N(b)\delta\right\}.
\end{array}
\]
Thus,
\[
\omega(\delta)=\inf\limits_{b>0}\left\{\frac{|\varphi(b)|}{|\psi(b)|}+N(b)\delta\right\}=l(\delta).
\]
Theorem is proved. $\square$

%
%
%
%
\begin{acknowledgements}
This project was supported by Simons Collaboration Grant N. 210363.
\end{acknowledgements}


\end{document}